\documentclass[eept,oneside,reqno]{amsart}
\usepackage{amsfonts}
\usepackage{amscd,amsmath}
\usepackage{url}
\usepackage{graphicx}
\usepackage{gastex}
\usepackage{longtable}
\usepackage{lscape}
\usepackage{tabularx}
\usepackage{multicol}
\usepackage{blindtext}
\usepackage{amsfonts}
\usepackage{amscd,amsmath}
\usepackage{url}
\usepackage{graphicx}
\usepackage{gastex}
\usepackage{longtable}
\usepackage{lscape}
\usepackage{tabularx}
\usepackage{multicol}
\usepackage{verbatim}
\usepackage{multirow}
\usepackage{graphicx}
\usepackage{epsfig,graphics,graphicx}
\usepackage{geometry}
\usepackage{verbatim}
\usepackage{multirow}
\usepackage{ragged2e}
\usepackage{tikz, tikz-cd}
\usepackage{amssymb}
\usepackage{setspace}
\newtheorem{theorem}{Theorem}[section]
\newtheorem{lemma}[theorem]{Lemma}

\usepackage[german,english]{babel}
\usepackage{multirow}
\usepackage{gensymb}

\newtheorem{corollary}[theorem]{Corollary}
\newtheorem{proposition}[theorem]{Proposition}
\newtheorem{definition}[theorem]{Definition}
\newtheorem{example}[theorem]{Example}

\theoremstyle{remark}
\newtheorem{remark}{Remark}[section]
\newtheorem*{rem*}{Remark}

\hbadness=\maxdimen

\numberwithin{equation}{section}
\usepackage{graphicx}
\usepackage{amsfonts}
\usepackage{tikz-cd}
\usepackage{mathtools}
\usepackage{amscd,amsmath}
\usepackage{amssymb}
\usepackage{xfrac}
\usepackage[all,cmtip]{xy}
\usepackage{amscd,amsmath}
\usepackage{amssymb}
\usepackage{fullpage}
\usepackage{setspace}
\theoremstyle{remark}
\numberwithin{equation}{section}

\begin{document}

\title[Multiplicative Lie algebra]{Isoclinism in multiplicative Lie algebras}
\author[Pandey]{Mani Shankar Pandey \vspace{.4cm}}
\address{The Institute of Mathematical Sciences, A CI of Homi Bhabha National Institute, C.I.T.Campus, Taramani, Chennai 600113, India}
\email{\tiny{manishankarpandey4@gmail.com}}.\\
\author[Upadhyay]{Sumit Kumar Upadhyay\vspace{.4cm}}
\address{Indian Institute of Information Technology, Allahabad, Prayagraj, 211015 India}
\email{\tiny{upadhyaysumit365@gmail.com}}.

\begin{abstract} 
The purpose of this paper is to introduce the notion of isoclinism and cover in a multiplicative Lie algebra which may be helpful to describe all multiplicative Lie algebra structures on a group. Consequently, we give the existence of the stem multiplicative Lie algebra.  We also give the necessary and sufficient condition for the existence of stem cover of a multiplicative Lie algebra.

\end{abstract}
\maketitle

\section{Introduction}
In the study of commutator identities of a group, G. J. Ellis (see \cite{GJ}) introduced a new algebraic structure and termed it multiplicative Lie algebra which is a generalization of group as well as Lie algebra. It appears that the concept of multiplicative Lie algebra has its own intrinsic interest. 

To see more algebraic properties of a  multiplicative Lie algebra, in 1996 Point and  Wantiez \cite{FP} gave the concept of nilpotency for multiplicative Lie algebras in parallel to Lie algebras and proved many nilpotency results which are known for groups and Lie algebras. In 2021, Pandey et al. \cite{MRS} discussed a different notion of solvable and nilpotent multiplicative Lie algebra with the help of  multiplicative Lie center and Lie commutator to tackle the question "how far a multiplicative Lie algebra from being an improper multiplicative Lie algebra". So many algebraic properties of a  multiplicative Lie algebra have studied in \cite{AGNM, GNM, GNMA, GM} which are already known for groups as well as multiplicative Lie algebra.
 In 2018, Lal and Upadhyay \cite{RLS} introduced the notion of Schur multiplier for a multiplicative Lie algebra as the second cohomology of multiplicative Lie algebras and proved that for a finite multiplicative Lie algebra $\mathcal{G}$ with presentation $1\longrightarrow \mathcal{R}\longrightarrow \mathcal{F}\longrightarrow \mathcal{G}\longrightarrow 1$, the factor  $\frac{\mathcal{R} \cap [\mathcal{F},\mathcal{F}](\mathcal{F} \star \mathcal{F})}{[\mathcal{R},\mathcal{F}](\mathcal{R}\star \mathcal{F})}$ is universal and they termed it  as the Schur multiplier $\mathcal{\Tilde{M}}(\mathcal{G})$ of $\mathcal{G}$.
 On a group $(\mathcal{G},\cdot)$, we have two apparent multiplicative Lie algebra structures, the trivial one given by $x\star y=1$, and improper one is given by $x\star y=[x,\ y]$, for all $x,y \in \mathcal{G}$. So it is natural to classify the multiplicative Lie algebra structures on a group.  In this direction, Pandey and Upadhyay \cite{MS1} tried to classify all possible multiplicative Lie algebra structures for finite groups. Continuing the line of investigation Walls \cite{Walls} has classified the multiplicative Lie algebra structures for many groups.

 Classification of groups is one of the oldest problem in theory of groups. In a classification of $p$-groups in  1940 P. Hall   \cite{Hall} introduced the notion of isoclinism of groups. Consequently, he proved that each isoclinism class of a finite group contains group $\mathcal{G}$ of minimal order with the additional condition, $Z(\mathcal{G})\subseteq [\mathcal{G},\mathcal{G}]$ and called this group $\mathcal{G}$ as the stem group. The finite-dimensional Lie algebra analog to isoclinism was developed in 1994 by K. Moneyhun \cite{K}. Alike groups, she defined the stem Lie algebra such that its center is contained in its derived subalgebra. This proves that for finite-dimensional Lie algebras the concepts of isoclinism and isomorphism are identical.

As an application of the concept of isoclinism for multiplicative Lie algebras, we see that isomorphism implies isoclinism but the converse need not be true, We also see that trivial multiplicative Lie algebra and improper multiplicative Lie algebra on a group are isoclinic. So, the concept of isoclinism may be helpful to classify all the multiplicative Lie algebra structures on a group.  

\section{Isoclinism in multiplicative Lie algebra}
The main aim of this section is to introduce the concept of isoclinism in multiplicative Lie algebras.
\begin{definition}
Let $(G,\cdot,\star)$ be a multiplicative Lie algebra. 
\begin{enumerate}
\item We call the ideal $[G,G](G\star G)$ as the multiplicative commutator of $G$ and denote it by $^M[G,G]$. For $g, g' \in G$, the element $^M[g,g^\prime]= (g\star g^\prime)[g, g^\prime]$ is called   the multiplicative commutator  of $g$ and $g^\prime$
\item We call the ideal $Z(G)\cap LZ(G)$ as the multiplicative Lie center of $G$ and denote it by $\mathcal{Z}(G)$. 
\end{enumerate}
\end{definition}

\begin{proposition} \label{improper}
The multiplicative commutator is the smallest ideal of $G$ such that the factor $\frac{G}{^M[G,G]}$ is an abelian group with trivial multiplicative Lie algebra structure.
\end{proposition}
\begin{proof}
To prove the theorem, it is enough to show that for any ideal $I$ of $G$,  $\frac{G}{I}$ is an abelian group with trivial multiplicative Lie algebra structure if and only if $^M[G,G]\subseteq I$. 

Let $I$ be an ideal of $G$. Then
\begin{center}
$\frac{G}{I}$ is abelian  $\Leftrightarrow [G,G] \subseteq I$. Also, $\frac{G}{I}$ is a trivial  multiplicative Lie algebra $\Leftrightarrow G\star G \subseteq I$
\end{center}
Hence, $\frac{G}{I}$ is an abelian group with trivial multiplicative Lie algebra structure if and only if $^M[G,G]\subseteq I$.
\end{proof}
\begin{proposition}\label{quotient}
Let $A$ be a multiplicative Lie algebra, $\Tilde{A}$ be an abelian Lie ring and $\rho: A \longrightarrow \Tilde{A}$ be a multiplicative Lie algebra homomorphism. Then there exists a unique multiplicative Lie algebra homomorphism $\bar{\rho}:\frac{A}{^M[A,A]}\longrightarrow \Tilde{A}$ such that the following diagram 
\begin{center}
\begin{tikzcd}
  A \arrow{r}{\rho} \arrow{d}{\nu}
    & \Tilde{A}  \\
  \frac{A}{^M[A,A]} \arrow{ur} {\exists~ \bar{\rho}}
&\end{tikzcd}
\end{center}
is commutative, where $\nu : A \longrightarrow \frac{A}{^M[A,A]}$ is the quotient multiplicative Lie algebra homomorphism.
Moreover, there exists a one to one correspondence between $Hom(A, \Tilde{A})$ and $Hom\big(\frac{A}{^M[A,A]}, \Tilde{A}\big)$.
\end{proposition}
\begin{proof}
Let $a,b,c,d \in A$. Then $\rho([a,b])=1$ and $\rho(c\star d)=1$. This shows that $ ^M[A,A]\subseteq ker(\rho)$. Thus we have a multiplicative Lie algebra homomorphism  $\bar{\rho}:\frac{A}{^M[A,A]}\longrightarrow \Tilde{A}$ defined by $\bar{\rho}(g\ ^M[A,A])=\rho(g)$. Clearly $\bar{\rho}$ is unique  and $\bar{\rho}\circ \nu (g)=\rho(g)$, for all $g \in A$.

Now, define $\Phi : Hom(\frac{A}{^M[A,A]}, \Tilde{A}) \longrightarrow Hom(A, \Tilde{A})$ by $\Phi(\rho)=\rho \circ \nu$. 
Let $\rho_1, \ \rho_2 \in Hom(\frac{A}{^M[A,A]}, \Tilde{A}) $ and $\Phi(\rho_1)=\Phi(\rho_2)$. Then $\rho_1 \circ \nu=\rho_2 \circ \nu\Rightarrow \rho_1= \rho_2$. So $\phi$ is one-one. Similarly, if $\rho \in Hom(A, \Tilde{A})$, then from above there exists $\bar{\rho} \in Hom(\frac{A}{^M[A,A]}, \Tilde{A})$ such that $\bar{\rho}\circ \nu=\rho$. So $\Phi(\bar{\rho})=\rho$. Thus we have a one to one correspondence between $Hom(\frac{A}{^M[A,A]}, \Tilde{A})$ and $Hom(A, \Tilde{A})$.
\end{proof}
\begin{remark}
Every abelian group can be treated as a trivial multiplicative Lie algebra and every group homomorphism can be treated as a trivial multiplicative Lie algebra homomorphism.
\end{remark}
Consider the category $\mathcal{ML}$ of multiplicative Lie algebras and category $\mathcal{AB}$ of abelian groups. Then we have two functors $\mathcal{F}: \mathcal{ML}\longrightarrow \mathcal{AB}$ and $\mathcal{I}: \mathcal{AB}\longrightarrow \mathcal{ML}$ given by 
\begin{center}
$\mathcal{F}(A)=\frac{A}{^M[A,A]}$ and $\mathcal{I}(A)=A_t$, 
\end{center}
where $A_t$ denotes the group $A$ with trivial multiplicative Lie algebra structure. 
\begin{theorem}
The functor $\mathcal{F}$ is a left adjoint  to the functor $\mathcal{I}$.
\end{theorem}
\begin{proof}
Let $M \in \mathcal{ML}$ and $A \in \mathcal{AB}$. Define $\Phi_{M,A}: Hom(\mathcal{F}(M), A) \longrightarrow Hom(M, \mathcal{I}(A))$ as $\Phi_{M,A} (\tau) = \tau \circ \nu$, for all $\tau \in Hom(\mathcal{F}(M), A)$. Then by Proposition \ref{quotient}, it is clear that $\Phi_{M,A}$ is a bijective map.

Now, let $M,\ \Tilde{M} \in \mathcal{ML}$, $f \in Hom(\Tilde{M},M)$ and $A,\ \Tilde{A} \in  \mathcal{AB}$, $g \in Hom(A,\Tilde{A})$.  Then we have two maps $\Psi_{M,\Tilde{M}}:Hom(\mathcal{F}(M),A)\longrightarrow Hom(\mathcal{F}(\Tilde{M}),\Tilde{A})$ and $\Psi_{A,\Tilde{A}}:Hom(M,\mathcal{I}(A))\longrightarrow Hom(\Tilde{M},\mathcal{I}(\Tilde{A}))$ defined by
\begin{center}
$\Psi_{M,\Tilde{M}}(\tau)=g\circ \tau \circ \mathcal{F}(f)$ \\
$\Psi_{A, \Tilde{A}}(\eta)=g\circ \eta\circ f$,
\end{center}
where $\mathcal{F}(f): \mathcal{F}(\Tilde{M})  \longrightarrow \mathcal{F}(M)$ is given by $\mathcal{F}(f)(\Tilde{M} \ {^{M }[\Tilde{M},\Tilde{M}]}) = f(\Tilde{M})\ {^{M} [M,M]}$.

\textbf{Claim:} $\Phi_{\Tilde{M},\Tilde{A}}\circ \Psi_{M,\Tilde{M}}=\Psi_{A,\Tilde{A}}\circ \Phi_{M,A}$, that is, the following diagram is commutative:
\[\begin{tikzcd}
Hom(\mathcal{F}(M),A) \arrow{r}{\Phi_{M,A}} \arrow[swap]{d}{\Psi_{M,\Tilde{M}}} & Hom(M,\mathcal{I}(A)) \arrow{d}{\Psi_{A,\Tilde{A}}} \\
Hom(\mathcal{F}(\Tilde{M}),\Tilde{A}) \arrow{r}{\Phi_{\Tilde{M},\Tilde{A}}} & Hom(\Tilde{M},\mathcal{\Tilde{A}}),
\end{tikzcd}
\] where $\Phi_{\Tilde{M},\Tilde{A}}: Hom(\mathcal{F}(\Tilde{M}), \Tilde{A}) \longrightarrow Hom(\Tilde{M}, \mathcal{I}(\Tilde{A}))$ defined by $\Phi_{\Tilde{M},\Tilde{A}} (\tau') = \tau' \circ \nu'$, for all $\tau' \in Hom(\mathcal{F}(\Tilde{M}), \Tilde{A})$ and $\nu' : \Tilde{M} \longrightarrow \frac{\Tilde{M}}{^L[\Tilde{M},\Tilde{M}]}$ is the quotient multiplicative Lie algebra homomorphism.

For this, let $\tau \in Hom(\mathcal{F}(M),A)$. Then 
$\Phi_{\Tilde{M},\Tilde{A}}(\Psi_{M,\Tilde{M}}(\tau))=\Phi_{\Tilde{M},\Tilde{A}}(g \circ \tau \circ \mathcal{F}(f))= g\circ \tau \circ \mathcal{F}(f) \circ \nu' = g\circ \tau \circ \nu \circ f $ (since $\mathcal{F}(f) \circ \nu' = \nu \circ f$). So $\Phi_{\Tilde{M},\Tilde{A}}(\Psi_{M,\Tilde{M}}(\tau)) = \Psi_{A,\Tilde{A}}(\tau \circ \nu)=\Psi_{A,\Tilde{A}}(\Phi_{M,A} (\tau) )$. Thus $\Phi_{\Tilde{M},\Tilde{A}}\circ \Psi_{M,\Tilde{M}}=\Psi_{A,\Tilde{A}}\circ \Phi_{M,A}$.
This proves that $\mathcal{F}$ is left adjoint to $\mathcal{I}$ or $\mathcal{I}$ is right adjoint to $\mathcal{F}$.
\end{proof}

The following lemma shows that the value of $^M[g,g^\prime]$ depends only on the pair of cosets of $\mathcal{Z}(G)$ to which $g,g^\prime$ belong.

\begin{lemma}
The maps $\phi_c:\frac{G}{\mathcal{Z}(G)}\times \frac{G}{\mathcal{Z}(G)} \longrightarrow \ ^M[G,G]$ and $\phi_s:\frac{G}{\mathcal{Z}(G)}\times \frac{G}{\mathcal{Z}(G)} \longrightarrow \ ^M[G,G]$ defined by $\phi_c(\tilde{g},\tilde{g}^{\prime})=\ [g,g^{\prime}]$ and $\phi_s(\tilde{g},\tilde{g}^{\prime})=\ (g\star g^{\prime})$ are well defined, where $\tilde{g},\tilde{g}^{\prime}$ are the images of $g,g^{\prime}$ in $\frac{G}{\mathcal{Z}(G)}$, respectively.
\end{lemma}
\begin{proof}
Let $\tilde{g},\tilde{g}^{\prime},\tilde{h},\tilde{h}^{\prime}\in \frac{G}{\mathcal{Z}(G)}$ and
$(\tilde{g},\tilde{g}^{\prime})=(\tilde{h},\tilde{h}^{\prime}).\ \text{Then}\  gh^{-1},g^{\prime}h^{\prime-1}\in \mathcal{Z}(G).$ Since $gh^{-1}\in LZ(G)$, we have
$
(g^{\prime}\star gh^{-1})=1 \Rightarrow  (g^{\prime} \star g)^g(g^{\prime}\star h^{-1})=(g^{\prime} \star g)\ ^{gh^{-1}}(g^{\prime}\star h)^{-1}=1$
Now, since  $gh^{-1}\in Z(G),\ \text{we have}\ (g^{\prime}\star g)=(g^{\prime}\star h).$ Thus, $$g\star g^{\prime} = h\star g^{\prime}.$$

Similarly, since $(h\star g^{\prime}h^{\prime-1})=1$, we have $(h\star g^{\prime})=(h\star h^{\prime}).$ Therefore, $$g\star g^{\prime} = h\star h^{\prime}.$$

On the other hand, since $gh^{-1}, g^{\prime}h^{\prime-1}\in Z(G)$, we have $[g^\prime,gh^{-1}]=1$ and $[h,g^{\prime}h^{\prime -1}] = 1$. This implies that $$[g, g^{\prime} ]= [h, h^{\prime}].$$
 Thus, $\phi_c$ and $\phi_s$ are well defined. 
\end{proof}
Now, we give definition of isoclinic multiplicative Lie algebras.
\begin{definition}
Two  multiplicative Lie algebras $G$ and $H$ are said  to be isoclinic (written as $G \sim_{ml} H$) if  there exist multiplicative Lie algebra isomorphisms $\lambda:\frac{G}{\mathcal{Z}(G)} \longrightarrow \frac{H}{\mathcal{Z}(H)}$ and $\mu:\ ^M[G, G]\longrightarrow \ ^M[H, H]$ such that the image of $^M[a,b]$ under $\mu$ is compatible with the image of $a\mathcal{Z}(G)$ and $b\mathcal{Z}(G)$ under $\lambda$. In other words the following diagram

\[
\xymatrix{
^M[G, G]\ar[d]^{\mu} & \frac{G}{\mathcal{Z}(G)}\times \frac{G}{\mathcal{Z}(G)}\ar[l]^{\phi_{c}~~~~~}  \ar[r]_{~~~~\phi_{s}} \ar[d]_{\lambda\times \lambda} &^M[G, G]  \ar[d]_{\mu} \\
^M[H, H] & \frac{H}{\mathcal{Z}(H)}\times \frac{H}{\mathcal{Z}(H)}\ar[l]_{\psi_{c}~~~~~~} \ar[r]^{~~~~\psi_{s}} & ^M[H, H]
}
\]

is commutative. The pair $(\lambda,\mu)$ is called as isoclinism between the  multiplicative Lie algebras $G$ and $H$.
\end{definition}
If two groups $G$ and $H$ are isoclinic, then we write it as $G \sim_{gr} H$. Now we are giving some examples.
\begin{example} \label{Example 1}
\begin{enumerate}
\item Let $G$ be a multiplicative Lie algebra. Then $G\sim_{ml} G \times A$, for every abelian group $A$ with trivial multiplicative Lie algebra structure.
\item Consider the group $A = \langle x, y, z \mid x^2 = y^2 = z^2 =1, xy =yx, yz =zy, xz = zx\rangle$ with multiplicative Lie algebra structure defined by $x\star y = x, x\star z = 1, y\star z = 1$. Then $\mathcal{Z}(A) = \{1,z\}$ and $^M[A, A] = \{1, x\}$. Now, consider the Klein four group $V_4 = \langle a, b \mid a^2 = b^2 =1, ab =ba\rangle$ with multiplicative Lie algebra structure defined by $a\star b = a$. Then  $\mathcal{Z}(V_4) = \{1\}$ and $^M[V_4, V_4] = \{1, a\}$. Then it is easy to see that $A \ncong V_4$ but $A \sim_{ml} V_4$.
\item Consider the multiplicative Lie algebras $(V_4,\cdot,\star)$ (as given in (2)) and $\mathbb{Z}_4$. Then $^M[V_4,V_4]=\{1,a\}$ and $^M[\mathbb{Z}_4,\mathbb{Z}_4]=\{1\}$. Therefore, $V_4 \nsim_{ml} \mathbb{Z}_4$ but $V_4 \sim_{gr} \mathbb{Z}_4$ $($for group isoclinism see \cite{Hall}$).$
\item The quaternian group $\mathcal{Q}_2=\langle x,y \ |\ x^2=y^2;\ y^4=1;\ yxy=x\rangle $ and the dihedral group $\mathcal{D}_4 =\langle x,y\ | \ x^2=y^4=1;\ xyx=y^3\rangle $ are isoclinic as group. It is well known from \cite{MS2} that $\mathcal{Q}_2$ is Lie simple. Consider $\mathcal{Q}_2$ with improper multiplicative Lie algebra structure and  $\mathcal{D}_4$ with $\star$ given by $x\star y=y$. Then $\mathcal{Q}_2 \nsim_{ml} \mathcal{D}_4$ but $\mathcal{Q}_2 \sim_{gr} \mathcal{D}_4$ .
\end{enumerate}
\end{example}
\begin{remark}
From the above examples, it can be seen that there is no direct relation between isoclinism of groups and isoclinism of multiplicative Lie algebras. But in the case of trivial and improper multiplicative Lie algebra structure, the notion of multiplicative Lie algebra isoclinism is same as that of groups. 
\end{remark}
\begin{remark}
For any group $G$, $G_t \sim_{ml} G_{I}$, where $G_t$ denotes  $G$ with trivial multiplicative Lie algebra  structure and $G_I$ denotes  $G$ with improper multiplicative Lie algebra  structure. 
\end{remark}
\begin{remark}
For perfect groups as well as perfect multiplicative Lie algebras  $G$ and $H$, $G \sim_{ml} H \Leftrightarrow G \cong H$ (since $^M[G,G] = G$ and $^M[H,H] = H$).
\end{remark}

\begin{lemma} \label{lem 1}
Let $G$ and $H$ be two multiplicative Lie algebras. Then
\begin{enumerate}
\item  $G$ is  a trivial Lie ring if and only if $G\sim_{ml} \{1\}.$ In particular  $G\sim_{ml} G\times \frac{G}{^M[G,G]}.$
\item $G\sim_{ml} H \Longleftrightarrow  G\times \cdots \times G\sim_{ml} H\times \cdots \times  H $.
\item If $H$ is a subalgebra of $G$, then $H\sim_{ml} H\mathcal{Z}(G)$. Moreover, $H\sim_{ml} G\Longleftrightarrow  H\mathcal{Z}(G)=G.$ 
\end{enumerate}
\end{lemma}
\begin{proof} The proof of (1) and (2) are obvious. We give the proof of (3).

   Since $\mathcal{Z}(G)$ is an ideal, $H\mathcal{Z}(G)$ is a subalgebra of $G$ and $\mathcal{Z}(H\mathcal{Z}(G)) = \mathcal{Z}(H)\mathcal{Z}(G)$. By Isomorphism theorem, $\frac{H\mathcal{Z}(G)}{\mathcal{Z}(H\mathcal{Z}(G))} \cong \frac{H\mathcal{Z}(H)\mathcal{Z}(G)}{\mathcal{Z}(H)\mathcal{Z}(G)} \cong \frac{H}{\mathcal{Z}(H)}.$ Thus, the map $\lambda:\frac{H}{\mathcal{Z}(H)}\longrightarrow \frac{H\mathcal{Z}(G)}{\mathcal{Z}(H\mathcal{Z}(G))}$  defined by $\lambda(h \mathcal{Z}(H))=h \mathcal{Z}(H\mathcal{Z}(G))$ is an isomorphism.
    
Now, let $h, h'\in H$ and $g, g'\in \mathcal{Z}(G)$. Then $(hg)\star (h'g') = h\star h'$ and $[(hg), (h'g')] = [h, h']$. Therefore, $^M[H\mathcal{Z}(G),H\mathcal{Z}(G)] \cong ^M[H,H]$ and the $\mu : ^M[H, H] \to ^M[H\mathcal{Z}(G),H\mathcal{Z}(G)]$ defined by $\mu([h, h']) = [h, h']$ and $\mu(h\star h') = h\star h'$ is an isomorphism. Now it is easy to see that the following diagram

\[
\xymatrix{
^M[H, H]\ar[d]^{\mu} & \frac{H}{\mathcal{Z}(H)}\times \frac{H}{\mathcal{Z}(H)}\ar[l]^{\phi_{c}~~~~~}  \ar[r]_{~~~~\phi_{s}} \ar[d]_{\lambda\times \lambda} &^M[H, H]  \ar[d]_{\mu} \\
^M[H\mathcal{Z}(G),H\mathcal{Z}(G)] & \frac{H\mathcal{Z}(G)}{\mathcal{Z}(H\mathcal{Z}(G))}\times \frac{H\mathcal{Z}(G)}{\mathcal{Z}(H\mathcal{Z}(G))}\ar[l]_{\psi_{c}~~~~~~} \ar[r]^{~~~~\psi_{s}} & ^M[H\mathcal{Z}(G),H\mathcal{Z}(G)]
}
\]

   is commutative. The next part of the proof follows directly from the case of groups (see \cite{MJ}  and \cite{Hall}).
\end{proof}
\begin{theorem}\label{theorem 1}
Let $G$ be a multiplicative Lie algebra and $I$ be an ideal of $G$. Then $\frac{G}{I}\sim_{ml}\frac{G}{I\cap \ ^M[G,G]}$. In particular, if $G$ is finite then $I\cap\ ^M[G,G]=\{1\}$ if and only if $G\sim_{ml} \frac{G}{I}$.
\end{theorem}
\begin{proof}
Suppose $\mathcal{Z}(\frac{G}{I}) = \frac{K}{I}$, where $K$ is an ideal of $G$. Let $aI \in \frac{K}{I}$. Then, $aI \star bI = I$ and $[aI, dI] = I$, that is, $a\star b \in I$ and $[a, b] \in I$, for every $b, d\in G$. Hence, $(G\star K)[G, K] \leq I$. Since $(G\star K)[G, K] \leq (G\star G)[G, G]$, $(G\star K)[G, K] \leq I\cap\ ^M[G,G]$. Hence, we have $\mathcal{Z}(\frac{G}{I}) = \frac{K}{I} \Leftrightarrow \mathcal{Z}(\frac{G}{I\cap ^M[G,G]}) = \frac{K}{I\cap ^M[G,G]}$. Therefore, $\frac{G}{I}/\mathcal{Z}(\frac{G}{I})\cong \frac{G}{I\cap \ ^M[G,G]}/\mathcal{Z}(\frac{G}{I\cap \ ^M[G,G]})$.

Since $[\frac{G}{I},\frac{G}{I}] \cong [\frac{G}{I\cap \ [G,G]},\frac{G}{I\cap \ [G,G]}]$ and $\frac{G}{I}\star \frac{G}{I} \cong \frac{G}{I\cap \ [G,G]}\star \frac{G}{I\cap \ [G,G]}$, $^M[\frac{G}{I},\frac{G}{I}] \cong \frac{^M[G,G]}{I\cap \ ^M[G,G]} =\ ^M[\frac{G}{I\cap \ ^M[G,G]},\frac{G}{I\cap \ ^M[G,G]}]$.   So, $\frac{G}{I}\sim_{ml}\frac{G}{I\cap \ ^M[G,G]}$.

If $I\cap\ ^M[G,G]=\{1\}$, then it is clear that $\frac{G}{I}\sim_{ml}G$. Conversely, suppose that $\frac{G}{I}\sim_{ml}G$.
 So, $\mid ^M[G,G] \mid =\ \mid \frac{^M[G,G]I}{I} \mid = \mid \frac{^M[G,G]}{I\cap \ ^M[G,G]} \mid \ \implies I \cap \ ^M[G,G]=\{1\}$. This proves the theorem.
\end{proof}
\begin{corollary} \label{cor 1}
Let $G$ and $H$ be two finite multiplicative Lie algebras and $\alpha$ be a homomorphism from $G$ into $H$, then
 $\alpha$ induces an isoclinism between $G$ and $H$ if and only if ker$(\alpha)\ \cap \ ^M[G,G] =\{1\}$ and
$H=$Image$(\alpha)\mathcal{Z}(H)$.
\end{corollary}
\begin{proof}
    This follows from Lemma \ref{lem 1} and Theorem \ref{theorem 1}.
\end{proof}
\begin{theorem}\label{theorem}
Let $G$ and $H$ be two multiplicative Lie algebras. Then $G$ and $H$ are isoclinic if and only
if there exists a multiplicative Lie algebra $\mathcal{K}$ containing ideals $\mathcal{Z}^G$ and  $\mathcal{Z}^H$ such that the following condition hold:
\begin{enumerate}
    \item $G\cong \frac{\mathcal{K}}{\mathcal{Z}^H} $
    \item $H\cong \frac{\mathcal{K}}{\mathcal{Z}^G} $
    \item $G\sim_{ml} \mathcal{K}\sim_{ml} H.$
\end{enumerate}
\end{theorem}
\begin{proof}
The result is obvious in one way.  We prove the converse part of the theorem. Let  $G$ and $H$ be isoclinic multiplicative Lie algebras and $(\lambda,\gamma)$ be an isoclinism between them. Let us consider a subalgebra $\mathcal{K}$ of the  multiplicative Lie algebra $G\times H$, $\mathcal{K}= \{(g,h) \  | \ \lambda(g\mathcal{Z}(G)=h\mathcal{Z}(H))\}$. Let $\mathcal{Z}^G=\{(g,1)\ |\ g\in \mathcal{Z}(G))\}$ and $\mathcal{Z}^H=\{(1,h)\ |\ h\in \mathcal{Z}(H))\}$. Then it is easy to verify that $\mathcal{Z}^G$ and $\mathcal{Z}^H$ are ideals of $\mathcal{K}$. Define the maps $\phi_{G}$ and $\phi_{H}$ as follows,
$$\phi_{G}:\mathcal{K}\longrightarrow G\ \text{and} \ \phi_{H}:\mathcal{K}\longrightarrow H$$ 
$$(g,h)\mapsto g\ \text {and}\ (g,h)\mapsto h$$
Then $\phi_G$ and $\phi_H$ are  surjective multiplicative homomorphisms with respective kernels $\mathcal{Z}^H$ and $\mathcal{Z}^G$. Therefore part $(1)$, and $(2)$ of the theorem has been done. It is easy to see  that the Lie abeleanizer $^M[\mathcal{K},\mathcal{K}]=\big\langle \big((g_1\star g_2)[g_3, g_4], \mu((g_1\star g_2)[g_3, g_4])\big)\ |\ g_i\in G, \ 1\leq i\leq 4\big\rangle$,   $^M[\mathcal{K},\mathcal{K}]\cap \mathcal{Z}^G=\{1\}=\ ^M[\mathcal{K},\mathcal{K}]\cap \mathcal{Z}^H$. Therefore by Theorem \ref{theorem 1}, $\frac{K}{\mathcal{Z}^H}\equiv \mathcal{K}\equiv \frac{K}{\mathcal{Z}^G}.$ This proves our assertion.
\end{proof}
The following theorem presents a relation between two isoclinic multiplicative Lie algebras. It asserts that any two isoclinic multiplicative Lie algebras  $G$ and $H$  have a common isoclinic descendant $\mathcal{K}$. 
\begin{theorem}\label{theorem 3}
Let $G$ and $H$ be two multiplicative Lie algebras. Then $G$ and $H$ are isoclinic if and only if there exists a multiplicative Lie algebra $\mathcal{\tilde{K}}$ containing ideals $\mathcal{\tilde{K}}_{{G}}$ and  $\mathcal{\tilde{K}}_{{H}}$ such that the following condition hold:
\begin{enumerate}
    \item $G\cong \mathcal{\tilde{K}}_G\sim_{ml}  \mathcal{\tilde{K}}$
    \item $H\cong \mathcal{\tilde{K}}_H \sim_{ml}  \mathcal{\tilde{K}}$
\end{enumerate}
\end{theorem}
\begin{proof} One way is obvious. Now, we prove the other way. Suppose $G\sim_{ml} H$. Then by Theorem \ref{theorem}, there exists a multiplicative Lie algebra
	 $\mathcal{K}$ containing ideals $\mathcal{Z}^G$ and  $\mathcal{Z}^H$ such that the following condition hold:
\begin{enumerate}
    \item $G\cong \frac{\mathcal{K}}{\mathcal{Z}^H} $
    \item $H\cong \frac{\mathcal{K}}{\mathcal{Z}^G} $
    \item $G\sim_{ml} \mathcal{K}\sim_{ml} H.$\end{enumerate}  Let $\mathcal{T}=\frac{\mathcal{K}}{\mathcal{Z}^H}\times \frac{\mathcal{K}}{^M[\mathcal{K},\mathcal{K}]}$. Then it can be verified that 
$$\mathcal{L}=\{((g,1)\mathcal{Z}^H,(g,1)^M[\mathcal{K},\mathcal{K}])\ |\ g\in \mathcal{Z}(G)\}$$ 
is an ideal of $\mathcal{T}$. Since $^M[\frac{\mathcal{K}}{\mathcal{Z}^H}, \frac{\mathcal{K}}{\mathcal{Z}^H}]\cong ^M[K, K]\cong ^M[\frac{\mathcal{K}}{\mathcal{Z}^G}, \frac{\mathcal{K}}{\mathcal{Z}^G}]$ and $^M[\mathcal{K},\mathcal{K}]\cap \mathcal{Z}^G=\{1\}=\ ^M[\mathcal{K},\mathcal{K}]\cap \mathcal{Z}^H$, $^M[\mathcal{T},\mathcal{T}]\cap \mathcal{L}=\{1\}$. Since for every $g\in G$ we have $h\in H$ such that $\lambda(g)\mathcal{Z}(G)=h\mathcal{Z}(H)$, we have multiplicative Lie algebra homomorphisms $\gamma_G : G\longrightarrow \frac{\mathcal{T}}{\mathcal{L}}\ \text{and}\ \gamma_H:H\longrightarrow \frac{\mathcal{T}}{\mathcal{L}}$ 
defined by  
$$\gamma_G{(g)}=\big((g,h)\mathcal{Z}^H,^M[\mathcal{K},\mathcal{K}]\big)\mathcal{L}\ \text{and} \ \gamma_H(h)=\big((g,h)\mathcal{Z}^H,(g,h)^M[\mathcal{K},\mathcal{K}]\big)\mathcal{L}.$$ 

Suppose $\gamma_G{(g)} = \mathcal{L}$. Then there exist $g'\in \mathcal{Z}(G)$ and $h\in H$ with $\lambda(g\mathcal{Z}(G))= h\mathcal{Z}(H)$ such that $(g'g^{-1}, h^{-1})\in \mathcal{Z}^H$ and $(g', 1) \in ^M[\mathcal{K},\mathcal{K}]$. Hence, $g'= g$, that is, $\gamma_G$ is injective. Similarly, $\gamma_H$ is also injective.

 Now, setting $\mathcal{\tilde{K}}=\frac{\mathcal{T}}{\mathcal{L}}$, $\mathcal{\tilde{K}}_G= \gamma_{G}(G)$ and $\mathcal{\tilde{K}}_H= \gamma_{H}(H)$. Clearly, $\mathcal{\tilde{K}}_G\cong G$ and $\mathcal{\tilde{K}}_H\cong H$.
 Now we claim that $\mathcal{\tilde{K}}=\mathcal{\tilde{K}}_G \mathcal{Z}(\tilde{K})$ and $\mathcal{\tilde{K}}=\mathcal{\tilde{K}}_H \mathcal{Z}(\tilde{K})$.

Let $(g_1,h_1),(g_2,h_2)\in \mathcal{K}$. Then $\big((g_1,h_1)\mathcal{Z}^H,(g_2,h_2)^M[\mathcal{K},\mathcal{K}]\big)\mathcal{L}\in \mathcal{\tilde{K}}$ and
\begin{multline*}
    \big((g_1,h_1)\mathcal{Z}^H,(g_2,h_2)^M[\mathcal{K},\mathcal{K}]\big)\mathcal{L}\\
=\bigg( \big((g_1,h_1)\mathcal{Z}^H,^M[\mathcal{K},\mathcal{K}] \big)\mathcal{L}\bigg)\bigg(\big( \mathcal{Z}^H, (g_2,h_2)^M[\mathcal{K},\mathcal{K}] \big)\mathcal{L}\bigg)\in \mathcal{\tilde{K}}_G \mathcal{Z}(\mathcal{\tilde{K}}),
\end{multline*}
Also
\begin{align*}
    ((g_1,h_1)\mathcal{Z}^H,(g_2,h_2)^M[\mathcal{K},\mathcal{K}])\mathcal{L}\\&=\bigg(((g_1,h_1)\mathcal{Z}^H,(g_1,h_1)^M[\mathcal{K},\mathcal{K}])\mathcal{L}\bigg)\bigg((\mathcal{Z}^H,(g_1^{-1}g_2,h_1^{-1}h_2)^M[\mathcal{K},\mathcal{K}])\mathcal{L}\bigg)\\ &\in \mathcal{\tilde{K}}_H \mathcal{Z}(\mathcal{\tilde{K}}).
\end{align*}
 Therefore, the theorem directly follows from Lemma \ref{lem 1} (3). 
%
\end{proof}
\begin{lemma}\label{lem 2}
Let $G$ and $H$ be two isoclinic multiplicative Lie algebras via isoclinism $(\lambda,\mu)$ and $g\in \ ^M[G,G],\ g^{\prime}\in G$. Then
\begin{enumerate}
\item $\lambda(g\mathcal{Z}(G))=\mu(g)\mathcal{Z}(H)$;
\item $\mu(^M[g,g^{\prime}])=\ ^M[\mu(g),h],$ for  $\lambda(g^{\prime}\mathcal{Z}(G)) = h\mathcal{Z}(H)$.  
\end{enumerate}
\end{lemma}
\begin{proof}
Let $g_1, g_2 \in G$ with $\lambda(g_1\mathcal{Z}(G)) = h_1\mathcal{Z}(H)$ and $\lambda(g_2\mathcal{Z}(G)) = h_2\mathcal{Z}(H)$. Then $\lambda((g_1 \star g_2)\mathcal{Z}(G)) = (h_1 \star h_2)\mathcal{Z}(H)$. Also, by the commutativity of the diagram, we have $\mu (g_1 \star g_2) = (h_1 \star h_2)$ . Therefore $\lambda((g_1 \star g_2)\mathcal{Z}(G)) = \mu (g_1 \star g_2)\mathcal{Z}(H)$.

Similarly, let $g_3, g_4 \in G$ with $\lambda(g_3\mathcal{Z}(G)) = h_3\mathcal{Z}(H)$ and $\lambda(g_4\mathcal{Z}(G)) = h_4\mathcal{Z}(H)$. Then $\lambda([g_3, g_4]\mathcal{Z}(G)) = ([h_3, h_4]\mathcal{Z}(H)$. Also, by the commutativity of the diagram, we have $\mu ([g_3, g_4]) = [h_3, h_4]$ . Therefore $\lambda([g_3, g_4]\mathcal{Z}(G)) = \mu ([g_3, g_4])\mathcal{Z}(H)$.

Thus, $\lambda((g_1 \star g_2)[g_3, g_4]\mathcal{Z}(G)) = \mu ((g_1 \star g_2)[g_3, g_4])\mathcal{Z}(H)$, for any $g_1, g_2, g_3, g_4 \in G$. Hence, $\lambda(g\mathcal{Z}(G))=\mu(g)\mathcal{Z}(H)$, for every $g\in \ ^M[G,G]$

%
Now, we prove the second part.  By part $(1)$, we have $\lambda(g\mathcal{Z}(G))=\mu(g)\mathcal{Z}(H)$. Thus by the commutativity of diagram, $\mu (g\star g')= \mu(g)\star h$ and $\mu ([g, g'])= [\mu(g), h]$. So, $\mu(^M[g,g^{\prime}])=\ ^M[\mu(g),h]$.
\end{proof}
Isoclinism partitions the family of all  multiplicative Lie algebras into disjoint classes. Any property depending on a variable multiplicative Lie algebra and which is the same for any two member of the same family will be called a class invariant.  Now we define the concept of stem multiplicative Lie algebra.
\begin{definition}
A multiplicative Lie algebra $G$ is said to be stem multiplicative Lie algebra in its equivalence class if $\mathcal{Z}(G)\subseteq \ ^M[G, G].$
\end{definition}
\begin{theorem} 
Let $G$ be a multiplicative Lie algebra. Then there exists a setm multiplicative Lie algebra $H$ such that $G$ and $H$ are isoclinic.
 \end{theorem}
\begin{proof}
Suppose $G$ is generated by set $\{g_i\ | \ i\in I \}$, where $I$ be an index set. Let $F$ be a trivial free Lie ring with the basis $\{f_i\ |\ i\in I\}$. Then $G\sim_{ml} G\times F$. Let $K$ be a subalgebra of $G\times F$ generated by the set $\{(g_i,f_i)\ | \ i\in I\}$. Then $KF$ coincides with $G\times F$. Since $F\subseteq \mathcal{Z}(KF)$, by Lemma \ref{lem 1}, $K\sim_{ml} G\times F \sim_{ml} G.$ Let $\pi$ be the projection map from $K$ to $F $ by mapping to the second coordinate and $\nu$ be the natural epimorphism from $K$ to $\frac{K}{^M[K,K]}$. Then there exits a  unique homomorphism $\phi$ from $\frac{K}{^M[K,K]}$ to $F$ such that the following diagram  
\begin{center}
\begin{tikzcd}
  K \arrow{r}{\pi} \arrow{d}{\nu}
    & F  \\
  \frac{K}{^M[K,K]} \arrow{ur} {\exists~ \phi}
&\end{tikzcd}
\end{center}
is commutative. Now, since $F$ is free abelian group, there exists a homomorphism $\hat{\phi}:F\longrightarrow \frac{K}{^M[K,K]}$ such that $\phi \circ \hat{\phi}=I$. In fact, $F\cong \frac{K}{^M[K,K]}$. Now,
\begin{align*}
\frac{\mathcal{Z}(K)}{\mathcal{Z}(K)\cap \ ^M[K,K]}\cong \frac{\mathcal{Z}(K)\ ^M[K,K]}{^M[K,K]}
\end{align*}
Therefore $\frac{\mathcal{Z}(K)}{\mathcal{Z}(K)\cap \ ^M[K,K]}$ is isomorphic to a subgroup of the free abelian group $F$. Thus, either $\frac{\mathcal{Z}(K)}{\mathcal{Z}(K)\cap \ ^M[K,K]}$ or $\frac{\mathcal{Z}(K)}{\mathcal{Z}(K)\cap \ ^M[K,K]}$ is free. 

Suppose $\frac{\mathcal{Z}(K)}{\mathcal{Z}(K)\cap \ ^M[K,K]}$. Hence, the theorem.  Now, suppose $\frac{\mathcal{Z}(K)}{\mathcal{Z}(K)\cap \ ^M[K,K]}$ is free. Thus $\mathcal{Z}(K)$ is the direct product of $\mathcal{Z}(K)\cap \ ^M[K,K]$ with another group $H$. Since $H\cap\ ^M[K,K]=\{1\} $, we have $\frac{K}{H}\equiv K$. Now,
 \begin{align*}
\mathcal{Z}\bigg(\frac{K}{H}\bigg)=\frac{\mathcal{Z}(K)}{H}=\frac{(\mathcal{Z}(K)\cap \ ^M[K,K])H}{H}\\
\subseteq \frac{ ^M[K,K]H}{H}\subseteq\ ^M\bigg[\frac{K}{H},\frac{K}{H}\bigg]
\end{align*}
and this proves the theorem.
\end{proof}
\begin{theorem} 
A multiplicative Lie algebras $G$ is stem in a class if and only  its order is minimal in that class.
\end{theorem}
\begin{proof}
The proof is similar as the case of groups (see \cite{MJ}  and \cite{Hall}).
\end{proof}

\begin{lemma}
Let $G$ and $H$ be two stem multiplicative Lie algebras in the same equivalence class. Then $\mathcal{Z}(G)\cong \mathcal{Z}(H).$
\end{lemma}
\begin{proof}
Let $(\lambda,\mu)$ be the isoclinism between $G$ and $H$. By the definition of stem multiplicative Lie algebra $\mathcal{Z}(G)\subseteq \ ^M[G,G]$. Let $g\in \mathcal{Z}(G)$. Then from Lemma \ref{lem 2}, $\lambda(g\mathcal{Z}(G))= \mu(g)\mathcal{Z}(H) = \mathcal{Z}(H)$. Thus $\mu(g)\in \mathcal{Z}(H)$. Thus the restriction of $\mu$ on  $\mathcal{Z}(G)$ is the required isomorphism between  $\mathcal{Z}(G)$ and $\mathcal{Z}(H)$.
\end{proof} 
 
 \section{Covers of multiplicative Lie algebras}
 Unlike the case of groups, K. Moneyhun \cite{K} has  shown that all covers of a finite-dimensional Lie algebra are isomorphic. A. R. Salemkar et. al. \cite{M} have shown the existence and uniqueness up to isoclinism of covers of any arbitrary  Lie algebra with finite-dimensional Schur multiplier.  In this section, we introduce a notion of stem cover of a multiplicative Lie algebra and we find a necessary and sufficient condition for stem cover of a multiplicative Lie algebra. 

 \begin{definition}
Let  $G$ be a multiplicative Lie algebra. Then a $2$-fold extension 
\[ 
\begin{tikzcd} 
1 \arrow{r} & H\arrow{r} & K\arrow{r} &G \arrow{r} & 1
\end{tikzcd}
\]
of $G$ is said to be a cover of $G$ if the following condition holds:
 \begin{itemize}
     \item $H\subseteq \mathcal{Z}(K)$
     \item $H\cong \mathcal{\Tilde{M}}(G).$
 \end{itemize}
 Further it is said to be an stem cover of $G$ if $H\subseteq \mathcal{Z}(K)\cap\  ^M[K,K]$.
 \end{definition}
 \begin{proposition}\label{l}
 Let \begin{tikzcd} 
  1 \arrow{r} &R\arrow{r} & F\arrow{r}{\pi} &G \arrow{r} & 1
\end{tikzcd} 
be a free presentation of a multiplicative Lie algebra $G$ and let 
\begin{tikzcd}  
1\arrow{r} & L\arrow{r} & K\arrow{r}{\beta} &\tilde{G} \arrow{r} & 1
\end{tikzcd}
be a central extension of another multiplicative Lie algebra $G'$. Then for each homomorphism $\phi$ from $G$ into $G$, there exists a homomorphism $\hat{\psi}$ from $\frac{F}{^M[R,F]}$ into $K$ such that $\hat{\psi}(R/^M[R, F ])$ is contained in $L$  and the following diagram 
\begin{center}
\begin{tikzcd}
    1\arrow{r} & \frac{R}{^M[R,F]}\arrow{r}{}\arrow{d}{\bar{\psi}} & \frac{F}{^M[R,F]}\arrow{r}{\bar{\pi}}\arrow{d}{\hat{\psi}} & {G}\arrow{r}\arrow{d}{\phi} & 1 \\
    1\arrow{r} & L\arrow{r}{} & K\arrow{r}{\beta} & {G'}\arrow{r} & 1
\end{tikzcd}
\end{center}

is commutative, where $\bar{\psi}$ is the restriction of the map $\hat{\psi}$ on $\frac{R}{^M[R,F]}$.
 \end{proposition}
 \begin{proof}
 As $F$ is free and $\beta$ is surjective, there is a homomorphism $\psi$ from $F$ into $K$ such that the following
diagram
 \[ \begin{tikzcd}
F\arrow{r}{\pi} \arrow[swap]{d}{\psi} & G \arrow{d}{\phi} \\%
K\arrow{r}{\beta}& G'
\end{tikzcd}
\]
 is commutative, where $\pi$ is the natural epimorphism. Now our aim is to show that the map $\hat{\psi}$ is induced by $\psi$. Since kernel of $\pi$ is $R$ and $\phi\circ \pi(r)=1=\beta\circ \psi(r)$ for each $r\in R$, we have $\psi(R)\subseteq \text{ker}(\beta)= L$. Now let $r\in R$ and $f\in F$. Then $ \psi(^M[r,f])= \ ^M[\psi(r),\psi(f)]$. Since $\psi(r)\in L$ and $L\subseteq \mathcal{Z}(K)$, we have   $\psi(^M[r,f])=1.$ Hence $\psi$ induces a homomorphism $\hat{\psi}$ from $\frac{F}{^M[R,F]}$ to $K$ as we required. 
 \end{proof}
 \begin{theorem}\label{th}
 Let $G$ be a multiplicative Lie algebra  and let 
 \begin{tikzcd} 
1 \arrow{r} & R\arrow{r} & F\arrow{r} &G \arrow{r} & 1
\end{tikzcd} 
be a free presentation of $G$. Then the  extension 
\begin{tikzcd} 
1 \arrow{r} & H\arrow{r} & G^{\star}\arrow{r}{\gamma} &G \arrow{r} & 1
\end{tikzcd} 
is a stem cover of $G$ if and only if there exists an ideal $T$ of $F$ such that following condition holds
\begin{enumerate}
    \item $G^\star \cong \frac{F}{T}$ and $H\cong \frac{R}{T}$;
    \item $\frac{R}{[R,F]}=\mathcal{\Tilde{M}}(G)\frac{T}{[R,F]}$, where $\mathcal{\Tilde{M}}(G)$ is the Schur multiplier of $G$.
\end{enumerate}
 \end{theorem}
 \begin{proof}
Suppose
 \[
 \begin{tikzcd} 
1 \arrow{r} & H \arrow{r} & G^{\star} \arrow{r}{\gamma} & G \arrow{r} & 1  \hspace{3cm} (1)
\end{tikzcd}  
\]
is a stem cover of $G$. Then by Proposition \ref{l}, there exists a homomorphism   $\tilde{\psi}$ from $\frac{F}{^M[R,F]}$ into $G^\star$ such that $\gamma \circ \tilde{\psi}=\bar{\pi}$ and $\Tilde{\psi}(\frac{F}{^M[R,F]}) \subseteq H$. Let $t$ be a transversal (for details about transversal see Section $3$ of \cite{MS1}) of the stem cover $(2)$. Then for every $g^\star\in G^\star$, there exists $h\in H$ and $x\in G$ such that $g^\star=ht(x)$. Again from the commutative diagram of Proposition \ref{l}, it can be seen that there exists $f\ ^M[R,F]\in \frac{F}{^M[R,F]}$ such that $\Tilde{\psi}(f\ ^M[R,F])=t(x)$. Therefore from the above discussion, it is easy to conclude that $G^\star =H\text{Im}(\Tilde{\psi})$. Since $(2)$ is a stem cover of $G$,  $H\subseteq \ ^M[G^\star,G^\star]=\ ^M[H\text{Im}(\Tilde{\psi}),H\text{Im}(\Tilde{\psi})]=\ ^M[\text{Im}(\Tilde{\psi}),\text{Im}(\Tilde{\psi})]$. This gives that $\Tilde{\psi}$ is onto and $\Tilde{\psi}(\frac{R}{^M[R,F]})=H.$ Choose an ideal $T$ of $F$ in such a manner that $\ker(\Tilde{\psi})=\frac{T}{^M[R,F]}$ and hence we get $G^\star \cong \frac{F}{T}$ and $H\cong \frac{R}{T}$. Now we prove the second assertion:
\begin{align*}
     \Tilde{\psi}(\mathcal{\Tilde{M}}(G))=  \Tilde{\psi}\bigg(\frac{R\cap \ ^M[F,F]}{^M[R,F]}\bigg)\subseteq \Tilde{\psi}\bigg(\frac{R}{^M[R,F]}\bigg)\cap \Tilde{\psi}\bigg(\frac{ ^M[F,F]}{^M[R,F]}\bigg) \subseteq H.
\end{align*}
 For any $h\in H$ there exists $r\ ^M[R,F]$ in $\frac{R}{^M[R,F]}$ such that $\Tilde{\psi}(r\ ^M[R,F])=h$. Since 
 $$H\subseteq \ ^M [G^\star,G^\star] \ =\ \Tilde{\psi}\bigg(\ ^M \big [\frac{F}{^M[R,F]},\frac{F}{^M[R,F]}\big]\bigg) , $$
we have some $f\ ^M[R,F]\in \frac{R}{^M[R,F]}$ such that $\Tilde{\psi}(r\ ^M[R,F])=\Tilde{\psi}(f\ ^M[R,F])=h$ Then $(rf^{-1}\ ^M[R,F])\in\ $ker$(\Tilde{\psi})=\frac{T}{^M[R,F]}\subseteq \frac{R}{^M[R,F]}$ and this gives that $f\in R$. Finally, we have $\Tilde{\psi}(\mathcal{\Tilde{M}}(G))=H.$ Now since the Schur multiplier is finite, $\mathcal{\Tilde{M}}(G)\cap \frac{T}{^M[R,F]}=\{1\}.$ This proves the second assertion. 
 
 Conversely suppose there exists an ideal $T$ of $F$ such that $H=\frac{R}{T}$ and $G^\star=\frac{F}{T}$. Then $\frac{G^\star}{H}=G$ and $H\cong \mathcal{\Tilde{M}}(G)$. Also, 
 \begin{align*}
     H=\frac{R}{T}\subseteq \frac{^M[F,F]T}{T}=\ ^M[G^\star,G^\star]
 \end{align*}
 Further, 
 \begin{align*}
     H=\frac{R}{T}\subseteq \mathcal{Z}(\frac{F}{T})=\mathcal{Z}({G^\star})
 \end{align*}
 This proves that the
extension \begin{tikzcd} 
1\arrow{r} & H\arrow{r} & G^{\star}\arrow{r} &G \arrow{r} & 1
\end{tikzcd} 
is a stem cover of $G$.
\end{proof}

\textbf{Acknowledgement} We would like to thank Prof. S. Vishwanath for his continuous support and discussion. Also, the  third named author wish to thank the Institute of Mathematical Sciences, Chennai for the financial support received through the institute's postdoctoral program.

\end{document}